\documentclass[12pt]{amsart}

\usepackage{color}
\usepackage{amsthm}
\usepackage{multicol}
\usepackage{amsfonts}
\usepackage{graphicx}
\usepackage{fancyhdr}
\usepackage{curves}

\newtheorem{theorem}{Theorem}[section]

\newtheorem{proposition}[theorem]{Proposition}

\newtheorem{definition}[theorem]{Definition}

\newtheorem{example}[theorem]{Example}

\begin{document}

\title{The geometry and combinatorics of Springer fibers}
\author{Julianna Tymoczko}

\begin{abstract}
This survey paper describes Springer fibers, which are used in one of the earliest examples of a geometric representation.  We will compare and contrast them with Schubert varieties, another family of subvarieties of the flag variety that play an important role in representation theory and combinatorics, but whose geometry is in many respects simpler.  The end of the paper describes a way that Springer fibers and Schubert varieties are related, as well as open questions.
\end{abstract}

\maketitle

\section{Introduction}

The flag variety $G/B$ is a critical object at the heart of geometry, combinatorics, algebra.  The interplay of these fields makes all of them richer.  In this survey article, we describe two families of subvarieties of the flag variety in some depth: Springer fibers and Schubert varieties.  We will compare the two families as well as point towards how they are related on a structural level.

Schubert varieties are the more foundational objects.  The double cosets $BwB$ partition $G$ into a finite number of subsets that the elements of the Weyl group $W$ naturally index.  These cosets induce a paving $BwB/B$ of the flag variety into Schubert cells, and the closure of $BwB/B$ is the Schubert variety corresponding to the Weyl group element $w$.  
The combinatorics of $w$ is intimately related to the geometry of the associated  Schubert variety.  Moreover the Schubert varieties induce a cohomology basis in $H^*(G/B)$ that plays an important role in representation theory.

Springer first defined the varieties now known by his name, as the fixed flags under the infinitesimal action of a nilpotent element $X$ \cite{Spr78}.  More precisely, the Springer fiber $\mathcal{S}_X$ associated to $X$ is the collection of flags
\[\mathcal{S}_X = \{gB \in G/B : g^{-1}Xg \in \mathfrak{b}\}\]
Springer used them to construct an action of the Weyl group on the cohomology $H^*(\mathcal{S}_X)$ that has many remarkable properties, described further in Section \ref{subsection: Springer representations}.  For this reason, and to understand more deeply the geometry underlying the construction, many have constructed Springer's representation (or its dual) \cite{KazLus80, BorMac83, Hot81, Lus84, GarPro92}.  These constructions often relied on geometric properties of Springer fibers like purity \cite{Spa76, Spa77}, closed formulas for dimensions \cite{Shi80, Fre09}, or identification of the components \cite{Spa76}.  The geometry is interesting in its own right, as well, and is connected to interesting combinatorial and representation-theoretic properties of $W$.

In what follows, we will describe these properties and connections in more detail, culminating in Section \ref{subsection: Springer representations} with a result that the Betti numbers of the Springer variety $\mathcal{S}_X$ are the same as a particular union of Schubert varieties determined by $X$.  We end with open questions.

Throughout we include many examples, focusing on examples in the type $A$ case because it is the most important geometric case.  

The author was partially supported by NSF grant DMS-1362855 and gratefully acknowledges the organizers of the conference at Orsay for the opportunity to discuss this work.

\section{Definitions and basic examples of Schubert varieties and Springer fibers}

Let $G$ be any complex semisimple connected linear algebraic group and let $B$ be a Borel subgroup of $G$.  The flag variety is the quotient $G/B$.  It is the ambient variety for everything that we consider in this paper.

The most important geometric example is the type $A$ case, when $G=GL_n(\mathbb{C})$ and $B$ is the subgroup of upper-triangular matrices.  In this case, each flag $gB$ can be thought of as a nested collection of subspaces
\[gB  \hspace{0.2in} \longleftrightarrow \hspace{0.2in} V_1 \subseteq V_2 \subseteq V_3 \subseteq \cdots \subseteq V_{n-1} \subseteq V_n = \mathbb{C}^n\]
where $V_i$ is an $i$-dimensional complex subspace for each $i \in \{1,2,\ldots,n\}$.  This correspondence is given by the rule that the first $i$ columns of $g$ span $V_i$ for each $i \in \{1,2,\ldots,n\}$. 

Now let $X$ be an element in the Lie algebra $\mathfrak{g}$ of $G$.  In type $A$  we can take $X$ to be a linear operator $X: \mathbb{C}^n \rightarrow \mathbb{C}^n$.  We write $E_{ij}$ for the basis element in $\mathfrak{gl}_n$ that has one in the $(i,j)$ position and zero elsewhere.

The Springer(-Grothendieck) fiber over $X$ is defined to be the collection of flags
\[\mathcal{S}_X = \{ gB: g^{-1}Xg \in \mathfrak{b} \}\]
Geometrically in type $A$ this can be rephrased as the collection of flags fixed by $X$ in the following sense:
\[\mathcal{S}_X = \{V_{\bullet}: XV_i \subseteq V_i \textup{ for all } i\}\]

The geometry of the flag variety is concrete and relatively well-adapted to direct computations.  For instance one classical way to analyze the flag $gB$ is via traditional ``Gaussian elimination", or column reduction, from undergraduate linear algebra.  Using column reduction, we can confirm the following.

\begin{proposition} \label{proposition: Schubert cell representatives}
Each coset $gB$ has a unique representative of the following form:
\begin{itemize}
\item a permutation matrix
\item with zeros below or to the right of its ones
\item and nonzero entries only above {\em and} to the left of its ones.
\end{itemize}
\end{proposition}

For instance one set of these representatives follows:
\[\begin{pmatrix} * & * & 1 & 0 \\ 1 & 0 & 0 & 0 \\ 0 & * & 0 & 1 \\ 0 & 1 & 0 & 0 \end{pmatrix}\]
The entries marked $*$ can be nonzero, and in fact can be any arbitrary element of $\mathbb{C}$.  This process hints at why the set of representatives that share a fixed permutation $w$ form a cell isomorphic to $\mathbb{C}^{d_w}$ for some integer $d_w$.  The cells $\{C_w\}$ are called Schubert cells.  In general types, this decomposition is precisely the one given by the double cosets:
\[C_w = BwB/B\]

Springer fibers are more complicated but similar arguments tell us quite a bit about their geometry.  If we choose a basis of $\mathbb{C}^n$ well with respect to $X$ then column reduction gives free entries in {\em some} of the free positions of $C_w$.  

\begin{example}\label{example: Springer cells} Let $X=E_{13}+E_{24}$ and consider the set of flags in the Schubert cell
\[\begin{pmatrix} * & * & 1 & 0 \\ 1 & 0 & 0 & 0 \\ 0 & * & 0 & 1 \\ 0 & 1 & 0 & 0 \end{pmatrix}\]
We can compute which of these flags are in $\mathcal{S}_X$ by calculating 
\[\begin{pmatrix} 0 & 0 & 1 & 0 \\ 0 & 0 & 0 & 1 \\ 0 & 0 & 0 & 0 \\ 0 & 0  &0 & 0 \end{pmatrix} \begin{pmatrix} a & b & 1 & 0 \\ 1 & 0 & 0 & 0 \\ 0 & c & 0 & 1 \\ 0 & 1 & 0 & 0 \end{pmatrix} = 
\begin{pmatrix} 0 & c & 0 & 1 \\ 0 & 1 & 0 & 0 \\ 0 & 0 & 0 & 0 \\ 0 & 0 & 0 & 0 \end{pmatrix} \]
and then determining the conditions that guarantee that the $i^{th}$ column in the matrix on the right is in the span of the first $i$ columns of the matrix on the left.  In this case, the condition is that $a=c$.
(In general the conditions on the entries can be more complicated.)
\end{example}

This basic approach allows us to determine entire Springer fibers, as in the following examples.

\begin{example}\label{example: Springers}
The first three examples are more general; the last is very specific.
\begin{itemize}
\item {\bf Suppose $X=0$.}  \\
\indent The condition that $g^{-1}Xg \in \mathfrak{b}$ is always satisfied so the Springer fiber $\mathcal{S}_0 = G/B$.
\item {\bf Suppose that $X$ is regular semisimple.}  (In type $A$ this means that $X$ has $n$ distinct eigenvalues.) \\
\indent The condition that $g^{-1}Xg \in \mathfrak{b}$ is satisfied if and only if $gB$ is a Weyl flag.  (In type $A$ we note that the first column of $g$ must be an eigenvector of $X$ and similarly the $i^{th}$ column of $g$ must be an eigenvector of $X$ for each $i$.  Putting this together implies that $gB$ is a permutation flag.)  This means $\mathcal{S}_X = \{wB: w \in W\}$.
\item {\bf Suppose that $X$ regular nilpotent.}  (In type $A$ this means that $X$ has a single Jordan block.) \\
\indent The condition that $g^{-1}Xg \in \mathfrak{b}$ is now satisfied if and only if $gB$ is the identity flag $eB$.  (In type $A$ the first column of $g$ must be an eigenvector of $X$ and similarly the $i^{th}$ column of $g$ must be an eigenvector of $X^i$ for each $i$.  Since $X$ has a unique eigenvector and more generally $X^i$ has $i$ eigenvectors there is a unique flag $gB$ satisfying the constraints.)  This means that $\mathcal{S}_X = \{eB\}$.
\item {\bf Suppose that $X=E_{13}$ in $\mathfrak{gl}_3$.} \\
\indent Concrete arguments like the previous tell us that the flag $g^{-1}Xg \in \mathfrak{b}$ if and only if $g$ is in one of the cosets
\[\begin{pmatrix}1 & 0 & 0 \\ 0 & 1 & 0 \\ 0 & 0 & 1 \end{pmatrix},
\begin{pmatrix} a  & 1 & 0 \\ 1 & 0 & 0 \\ 0 & 0 & 1 \end{pmatrix},
\begin{pmatrix}1 & 0 & 0 \\ 0 & b & 1 \\ 0 & 1 & 0 \end{pmatrix}\]
where $a,b \in \mathbb{C}$.  Inspecting what happens to the flags as each of $a,b$ gets large, we can see the closure of each one-(complex-)dimensional cell is the single point $eB$.  This means the Springer fiber is
\[\mathcal{S}_X = \left\{ \begin{array}{c}\textup{two spheres joined at }
\textup{ a single common point}\end{array} \right\}
\]
\end{itemize}
\end{example}

{\bf Unless otherwise indicated, we assume henceforth that $X$ is nilpotent.}

\section{Comparing the geometry of Schubert varieties and Springer fibers}

In the previous section we used linear algebraic techniques to identify the flags in specific Schubert varieties and Springer fibers.  In each case we found that the variety was partitioned into cells that were indexed by permutation flags.  In fact the combinatorics of the permutation at the heart of each cell determines quite a bit of the geometry of the cell and its closure.  The key difference between Schubert varieties and Springer fibers is that all entries in the cells of Schubert varieties are free while entries in cells of Springer fibers have conditions on them.  In this section we describe the ramifications of this key difference, from its impact on the dimension of the cells and how they intersect, to the kinds of open questions we have about each.  We start by sketching the main properties of Schubert cells and Schubert varieties, and then proceed to describe the cells in Springer fibers.  The reader interested in learning more about the geometry of Schubert varieties is referred to texts like those of Brion \cite{Bri05}, Fulton \cite{Ful97}, or Billey and Lakshmibai \cite{BilLak00}.

\subsection{Geometry of Schubert cells in the flag variety}
Decomposing the flag variety into  Schubert cells $C_w=BwB/B$ is very natural from an algebraic perspective.  Amazingly, this decomposition is more than just a partition of the flag variety: the cells actually form a CW-decomposition of the flag variety.  In other words
\begin{itemize}
\item each Schubert cell $C_w$ is isomorphic to $\mathbb{C}^{\ell(w)}$ for some integer $\ell(w)$ (as evident from the construction above) and 
\item the closure $\overline{C_w}$ of each Schubert cell is a union of Schubert cells of smaller dimension (less evident, but since each $C_w$ is a $B$-orbit, its closure is a union of $B$-orbits, namely other Schubert cells $C_v$). 
\end{itemize} 
More amazingly yet, the statistics that guarantee the CW-decomposition are completely natural combinatorial statistics.  The dimension $\ell(w)$ is the {\em length} of the permutation $w$, namely the smallest number of simple reflections $s_i=(i,i+1)$ into which $w$ can be factored.  (Length is defined similarly for the other Weyl groups.)  The closure relation is given by 
\[\overline{C_w} = \bigcup_{v \leq w} C_v\]
where the partial order $v \leq w$ is the {\em Bruhat order}.  The combinatorial definition of $v \leq w$ in Bruhat order is that $v$ is a subword of $w$ when $w$ is written in terms of the simple reflections $s_i=(i,i+1)$.  For instance if we consider the permutation $w=3142$ in one-line notation (meaning that $w(i)$ is in the $i^{th}$ position) then we can factor $w=s_2s_3s_1$.  The Weyl group element  $s_2s_3s_2s_1$ has length $4$ and is greater than the following set of Weyl group elements:
\[\{ e, s_1, s_2, s_3, s_2s_3, s_2s_1, s_3s_2, s_3s_1, s_2s_3s_2, s_2s_3s_1, s_3s_2s_1\}\]

From a topological perspective, the most useful fact about CW-decompositions is that they induce cohomology bases.  In particular the closures $\{\overline{C_w}: w \in W\}$ induce a module basis for the cohomology $H^*(G/B, \mathbb{Q})$.  The closure $\overline{C_w}$ of the Schubert cell $C_w$ is called a Schubert variety and the classes they induce in $H^*(G/B, \mathbb{Q})$ are called Schubert classes.  They are one of the most important bases of the cohomology of the flag variety (and more generally partial flag varieties $G/P$ when $P$ is a parabolic subgroup) for reasons we will describe in Section \ref{subsection: Schubert calculus}.

Moreover the geometry of Schubert varieties can be used to determine the cohomology ring structure directly.  From a geometric perspective, the key point of intersection theory is that when varieties lie in appropriate relative positions, their intersection induces the product of the corresponding cohomology classes.  This basic principle is true but more complicated when the varieties have singularities, so understanding singularities is an important part of geometric calculations in cohomology.   

For this reason many researchers have studied the singularities of Schubert varieties.  What they discovered was that these singularities are also deeply entwined with the combinatorics of permutations, specifically patterns inside permutations.  Lakshmibai and Sandhya proved that the variety $\overline{C_w}$ is smooth if and only if the permutation $w$ avoids the patterns $3412$ and $4231$ \cite{LakSan90}.  The permutation $w$  avoids a pattern if, when $w$ is written in one-line notation, no four numbers have the same relative positions as the pattern.  (Both Ryan \cite{Rya87} and Wolper \cite{Wol89} characterized the singular Schubert varieties before Lakshmibai and Sandhya, and without this particular combinatorial formulation; all three results are independent.)

\begin{example}
The permutations
\[w = 624351 \textup{      and     } v = 324651\]
both avoid the pattern $3412$ because only $4$, $5$, or $6$ can be the fourth-largest number in a subset of $\{1,2,3,4,5,6\}$ and none of them can have three smaller numbers in the appropriate relative positions.  The permutation $v$ also avoids the pattern $4231$ and so $\overline{C_v}$ is smooth.  The permutation $w = {\bf 6 2} 4 3 {\bf 5} {\bf 1}$ does not, so $\overline{C_w}$ is singular.
\end{example}

Pattern avoidance is an important condition in computer science and enumerative combinatorics. It also turns out to be related to many other properties of Schubert varieties, for instance the components of the singular locus of a Schubert variety \cite{BilWar03, Cor01, KLR03, Man01}, or whether the Schubert variety $\overline{C_w}$ is Gorenstein (a geometric criterion that is not as strong as being smooth) \cite{WooYon06}. 
(See Woo and Yong's work for a unified presentation and extension of these results \cite{WooYon08}.)

\subsection{Geometry of cells in Springer fibers} \label{section: geometry of Springer cells}
Recall from Example \ref{example: Springer cells} that column reduction of flags in Springer fibers produced cells that satisfied the same conditions for Schubert cells  described in Proposition \ref{proposition: Schubert cell representatives} {\em except that not all entries needed to be free}.  This is because each of these cells in the Springer fiber is an intersection with a specific Schubert cell. Of course, intersections can be fantastically complicated depending on the conditions on the non-free entries.  In this section, we describe how this affects what is known about the geometry of the cells in Springer fibers.

\begin{figure}[h]
\begin{picture}(220,60)(0,-30)

\thicklines
\put(20,0){\color{red} \circle{100}}
\put(61,0){\color{green} \circle{100}}
\put(102,0){\color{blue} \circle{100}}
\put(143,0){\color{yellow} \circle{100}}
\put(184,0){\circle{100}}

\thinlines
\put(0,0){\color{red} {\curve(0,0, 20,-7, 40,0)}}
\put(41,0){\color{green} {\curve(0,0, 20,-7, 40,0)}}
\put(82,0){\color{blue} {\curve(0,0, 20,-7, 40,0)}}
\put(123,0){\color{yellow} {\curve(0,0, 20,-7, 40,0)}}
\put(164,0){{\curve(0,0, 20,-7, 40,0)}}
\curvesymbol{\phantom{\circle*{2}} \circle*{1}}
\put(-2,0){\color{red} {\curve[-6](0,0, 20,4, 40,0)}}
\put(39,0){\color{green} {\curve[-6](0,0, 20,4, 40,0)}}
\put(80,0){\color{blue} {\curve[-6](0,0, 20,4, 40,0)}}
\put(121,0){\color{yellow} {\curve[-6](0,0, 20,4, 40,0)}}
\put(162,0){{\curve[-6](0,0, 20,4, 40,0)}}

\put(0,0){\circle*{6}}
\put(40,0){\color{red} \circle*{4}}
\put(81,0){\color{green} \circle*{4}}
\put(122,0){\color{blue} \circle*{4}}
\put(163,0){\color{yellow} \circle*{4}}

\end{picture}
\caption{Paving a ``string of pearls" by affines}\label{figure: paving example}
\end{figure}
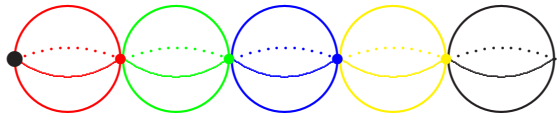

The first complication is that the intersections $\{C_w \cap \mathcal{S}_X\}$ do not form a CW-decomposition of the Springer fibers, unlike the Schubert cells for the flag variety.  Instead these intersections $\{C_w \cap \mathcal{S}_X\}$ form a topological decomposition called a {\em paving by affines}.  To be paved by {\em affines} means having a partition into cells that are each affine, namely isomorphic to $\mathbb{C}^{d_w}$ for some integer $d_w$.   (In other applications, it can be useful to study pavings by other geometric spaces.)  A {\em paving} restricts the closure conditions on each cell more loosely than in a CW-decomposition.  Cells form a paving if they can be ordered by an index set $\mathcal{I}$ so that the closure $\overline{C_i} \subseteq \bigcup_{j \leq i} C_j$ for each $i$ in the index set.  That containment---which would be an equality if the cells formed a CW-decomposition---changes the intuition around pavings by affines.  The classic example of a paving by affines that is not a CW-decomposition is a ``string of pearls", namely a collection of copies of $\mathbb{P}^1$ with the north pole of one glued to the south pole of the next.  The homologically natural way to partition this space is by pulling off the leftmost north pole, then the copy of $\mathbb{C}$ left in the leftmost $\mathbb{P}^1$, then the copy of $\mathbb{C}$ left in the next $\mathbb{P}^1$, and so on, as shown in Figure \ref{figure: paving example}.  This paving has one affine cell for each (co)homology basis class, with dimensions of the cells matching the degrees of the classes.  However it is not a CW-decomposition since the closure of most of the $\mathbb{C}$ cells contains a single point in the middle of another $\mathbb{C}$ cell. Historically Spaltenstein first described the irreducible components of each Springer fiber \cite{Spa76}; Shimomura refined his analysis to pave Springer fibers with affines \cite{Shi80}.  In both works, the dimension of the cells was naturally associated to certain combinatorial objects that we describe next.

The dimension $d_w$ of each cell can be determined combinatorially for Springer fibers but we only count part of the dimension $\ell(w)$ of the corresponding Schubert cell.  Another way of defining the length of a permutation $w$ is as the number of inversions of $w$.  If the permutation is written in one-line notation, an inversion is a pair of numbers with the bigger number to the left of the smaller number. For instance the permutation
\[ w = 513264\]
has $6$ inversions, one from the pair $3>2$, four from the numbers less than $5$, and one from the pair $6>4$.  (The definition of inversions can be extended to arbitrary Lie type using the combinatorics of roots.)    

\begin{figure}[h]
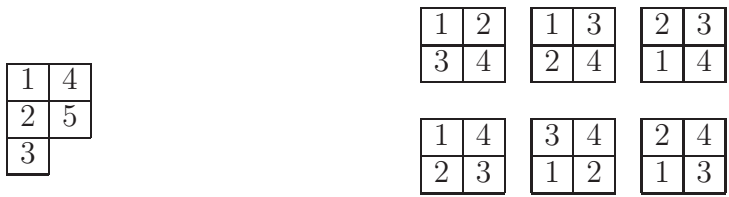

\[\begin{array}{cc}
\begin{array}{|c|c|} \cline{1-2} 1 & 4  \\ \cline{1-2} 2 & 5  \\ \cline{1-2} 3 & \multicolumn{1}{c}{} \\ \cline{1-1} \end{array}

&

\begin{array}{ccc}
\begin{array}{|c|c|}
\cline{1-2} 1 & 2 \\
\cline{1-2} 3 & 4 \\
\cline{1-2} \end{array}
&
\begin{array}{|c|c|}
\cline{1-2} 1 & 3 \\
\cline{1-2} 2 & 4 \\
\cline{1-2} \end{array}
&
\begin{array}{|c|c|}
\cline{1-2} 2 & 3 \\
\cline{1-2} 1 & 4 \\
\cline{1-2} \end{array}
\\
\hspace{1em} & & \\
\begin{array}{|c|c|}
\cline{1-2} 1 & 4 \\
\cline{1-2} 2 & 3 \\
\cline{1-2} \end{array}
&
\begin{array}{|c|c|}
\cline{1-2} 3 & 4 \\
\cline{1-2} 1 & 2 \\
\cline{1-2} \end{array}
&
\begin{array}{|c|c|}
\cline{1-2} 2 & 4 \\
\cline{1-2} 1 & 3 \\
\cline{1-2} \end{array}
 \\
\hspace{1em} & & \\
\end{array} \\
 
 \textup{     (a)   \hspace{0.5em} Young diagram example } \hspace{0.35in} &
 \textup{ (b)   \hspace{0.5em}   All row-strict Young tableaux } \\
 \textup{    } & 
 \textup{  of shape }(2,2) \\
\end{array}\]
\caption{Young tableaux}\label{figure: examples of Young tableaux}
\end{figure}

To count the dimension of the cells in Springer fibers, we need to incorporate both inversions and the structure of the matrix $X$ into the combinatorics.  To do this, we use a Young diagram, which is a top-aligned and left-aligned collection of $n$ boxes. (This is known as English notation rather than French notation, which traditionally uses bottom-aligned and left-aligned boxes.  To borrow Macdonald's quip \cite[page 2]{Mac95}, the reader accustomed to French notation should read this paper upside down in a mirror.)   In our case we take $\lambda(X)$ to be the Young diagram with the same number of rows as Jordan blocks in $X$ and the same number of boxes in each row as  the dimension of the corresponding Jordan block.  ``Top-aligned" means that we arrange the diagram so that rows decrease in size.  Figure \ref{figure: examples of Young tableaux} shows examples of Young diagrams, in this case filled with numbers.  We fill our Young diagrams bijectively with the integers $\{1,2,\ldots,n\}$, meaning that each integer appears exactly once.  A row-strict filling is one in which the integers increase from left-to-right in each row, without any conditions on columns.  

\begin{theorem} \label{theorem: all Springer cells}
 If $X$ is chosen carefully within its conjugacy class, the nonempty intersections $C_w \cap \mathcal{S}_X$ are bijective with the set of row-strict fillings of $\lambda(X)$.
\end{theorem}

For certain $X$ the cells $C_w \cap \mathcal{S}_X$ are in fact described geometrically by the fillings. Loosely speaking, reading the filling in a specific order gives a permutation, which tells in what order to add basis vectors to the flag (equivalently which permutation matrix appears in the coset representative).  A more precise argument can be found in \cite{Shi80} and \cite{Tym06} or \cite{Fre09}.  (Little is known about the intersections $C_w \cap \mathcal{S}_X$ for other $X$ except that they can be very complicated.  For instance in the analogous case of Peterson varieties, if $X$ is a particular lower-triangular matrix then the coordinate ring of the largest intersection is the quantum cohomology of the flag variety \cite{Kos96, Rie03}, while if $X$ is the opposite upper-triangular matrix then the intersections are affine cells \cite{Tym06}.)

Moreover the row-strict Young diagrams can be used to compute the dimensions of the corresponding intersection $C_w \cap \mathcal{S}_X$.  For each integer $i$ in the filling, count the integers $j$ that satisfy these conditions:
\begin{itemize}
\item $i>j$
\item $j$ is in any of the columns to the right of $i$ or is in the same column as $i$ but above $i$
\item either $j$ is the rightmost box in a row or the box to the right of $j$ is filled with a number $k$ that satisfies $i < k$
\end{itemize}
The first two conditions together are the conditions for $i$ and $j$ to form an inversion, if you read the entries in the Young diagram starting at the bottom of the left column, proceeding up each column, and then moving to the bottom of the next-right column.  The last condition gives a way to eliminate certain inversions---exactly ones that correspond to the dependent entries in the cell 
\cite{Tym06, Fre09}.  For example, the leftmost tableau in Figure \ref{figure: examples of Young tableaux} corresponds to a cell of dimension $4$.  The three tableaux on the top row of the right side of Figure \ref{figure: examples of Young tableaux} correspond to cells of dimension $2, 2, 1$ while the tableaux on the bottom row correspond to cells of dimension $1, 1, 0$.

Up to this point, the story of the cells in Springer fibers is analogous to Schubert cells in the flag variety,  but this is where the analogy breaks down.  Almost nothing is known about the geometry of the irreducible components of Springer fibers, much less the geometry of the closures of each cell $C_w \cap \mathcal{S}_X$.  Even relatively simple questions about the closure conditions are mysterious: for instance, we do not know which permutation flags lie in the boundary of the cell $C_w \cap \mathcal{S}_X$. Most of what is known about smoothness conditions consists of criteria determining when {\em all} components are smooth, which is of course a high standard to meet.  Indeed all components are smooth very rarely, including when $\lambda(X)$ has a hook shape \cite{Var79}, when $X$ has exactly two Jordan blocks \cite{Fun03}, and a small list of other cases \cite{FreMel10}.  The case when $X$ has exactly two Jordan blocks was particularly interesting to Khovanov because of connections to categorification in knot theory and quantum representations (see Section \ref{subsection: Khovanov Springers} and \cite{Kho04}).  Of course, if that little is known about which components are singular, it is unsurprising to learn that nothing is known about more refined conditions like being Gorenstein.  The exception is when $\lambda(X)$ has two columns, in which case the components are normal, Cohen-Macaulay, and rationally smooth \cite{PerSmi12} and smooth precisely when their Poincar\'{e} polynomials are symmetric \cite{FreMel11}.

Figure \ref{figure: comparison table} summarizes the discussion in this section.

\begin{figure}[h]
\begin{tabular}{ll}
  Flag variety $G/B$  & Springer fiber $\mathcal{S}_X$  \vspace{0.2em}\\
 \cline{1-2} & \vspace{-.5em}\\
  Partitioned into cells & Partitioned into affine pieces  \\
\hspace{0.25in} $C_w = BwB/B$ & \vspace{0.35em}\\
CW-decomposition into $C_w$ & Paved by affines $C_w \cap \mathcal{S}_X$ \vspace{0.35em} \\
$\{\overline{C_w}\}$ induce basis of $H^*(G/B)$ & $\{\overline{C_w \cap \mathcal{S}_X}\}$ induce basis of $H^*(\mathcal{S}_X)$ \vspace{0.15in} \\
 $\overline{C_w}= \bigcup_{v \leq w} C_w$  & ? \\
\hspace{-1em}  \begin{tabular}{l} smooth iff $w$ avoids 1324 and 2143 \\ 
\end{tabular} & \hspace{-1em} \begin{tabular}{l} \\ all components smooth: \\
$\bullet$ if $\lambda(X)$ has a hook shape \\
$\bullet$ if $X$ has two Jordan blocks 
\\ $\bullet$ for small list of shapes of $X$ \\ 
 \end{tabular} \vspace{0.15em}  \\
\hspace{-1em}  \begin{tabular}{l} \\ Gorenstein iff w avoids \\ 31524 and 24153 (+ conditions) 
\end{tabular} & ?\vspace{0.35em} \\

\end{tabular}
\caption{Comparing Schubert varieties and Springer fibers}\label{figure: comparison table}
\end{figure}

\section{Schubert varieties and Springer fibers in representation theory}

Combinatorics, geometry, and representation theory collide in the representations associated to Schubert varieties and Springer fibers.  Despite the similarities, the representations involve very different constructions---and groups!---using very different tools and approaches.  This section introduces three representations: the representation of $GL_n$ arising from the Schubert basis in the cohomology of the Grassmannian, the representation of $S_n$ on the cohomology of Springer fibers (or in general Lie type the representation of the Weyl group $W$), and more recent research into a quantum representation encoded by the components of certain Springer fibers.  Fulton's text is an excellent introduction to Schubert calculus \cite{Ful97} while Chriss and Ginzburg's is an expansive introduction to Springer theory \cite{ChrGin97}.

\subsection{Classical representation theory of Schubert varieties}\label{subsection: Schubert calculus}  
\begin{figure}[h]
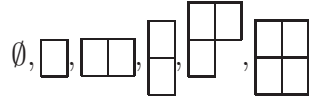

\[
\emptyset, 
\begin{array}{|c|}
\cline{1-1} \\
\cline{1-1} \end{array}, 
\begin{array}{|c|c|}
\cline{1-2} & \\
\hline \end{array}, 
\begin{array}{|c|}
\cline{1-1} \\
\cline{1-1} \\
\cline{1-1} \end{array}, 
\begin{array}{|c|c|}
\cline{1-2} & \\
\cline{1-2} & \multicolumn{1}{|c}{} \\
\cline{1-1} \multicolumn{2}{c}{}
 \end{array}, 
 \begin{array}{|c|c|}
\cline{1-2} & \\
\cline{1-2} & \\
\hline \end{array}
\]
\caption{The Young diagrams indexing Schubert classes in $H^*(G(2,4))$} \label{figure: Grassmannian Schubert basis}
\end{figure}

The most well-known and well-understood geometric representation involving Schubert varieties does not consider the Schubert classes within the cohomology of the flag variety.  Rather, it uses the image of the Schubert classes under the projection $G/B \rightarrow \hspace{-0.75em} \rightarrow G/P$ where $P$ is a maximal parabolic of type $A$.  In that case $G/P$ is a Grassmannian $G(k,n)$, namely the collection of $k$-dimensional subspaces of the vector space $\mathbb{C}^n$.   The surjection $G/B \rightarrow \hspace{-0.75em} \rightarrow G/P$ collapses many of the Schubert varieties in the full flag variety, in the sense that their images are not distinct. Instead of being indexed by the set of all permutations, the Schubert cells $\{C_{\lambda}\}$ in the Grassmannian are indexed by the Young diagrams $\lambda$ with at most $k$ columns and at most $n-k$ rows.  Figure \ref{figure: Grassmannian Schubert basis} shows the Young diagrams that index the Schubert cells for $G(2,4)$.  As this example demonstrates, the Grassmannian is generally much simpler than the full flag variety, since e.g. the variety $GL_4(\mathbb{C})/B$ has $24=4!$ Schubert classes instead of just $6$.

Now suppose $\sigma_{\lambda}$ represents the cohomology class of $\overline{C_{\lambda}}$ in $G(k,n)$.  We can write the product $\sigma_{\lambda} \cdot \sigma_{\mu}$ in the cohomology ring $H^*(G(k,n), \mathbb{C})$ in terms of its basis of Schubert classes as so:
\[\sigma_{\lambda} \cdot \sigma_{\mu} = \sum_{\nu} c^{\nu}_{\lambda, \mu} \sigma_{\nu}\]
It turns out that the coefficients $c^{\nu}_{\lambda, \mu}$ are precisely the same coefficients obtained in the tensor product decomposition
\[V_{\lambda} \cdot V_{\mu} = \sum_{\nu} c^{\nu}_{\lambda, \mu} V_{\nu}\]
where $V_{\lambda}$ are the irreducible representations of $GL_n(\mathbb{C})$.  

Moreover the combinatorics of Young diagrams determines many properties of these coefficients.  For instance if $\mu$ is not contained in $\nu$ then the coefficient $c^{\nu}_{\lambda, \mu} $ must be zero; similarly certain fillings of the Young diagrams $\lambda$ and $\mu$ determine the coefficients $c^{\nu}_{\lambda, \mu}$. (The fillings in Schubert calculus are not the same fillings that we used for Springer fibers.)  Some of these properties can be proven with very elegant geometry, while for others the only known proofs rely on technical results from the theory of symmetric functions.  Indeed the classical proof that the structure constants of $H^*(G(k,n))$ equal the tensor product multiplicities simply observes that both objects are determined by a recurrence relation counted by certain combinatorial objects, and then shows that the recurrence and the initial cases are the same in both contexts.  In other words, the proof is not terribly illuminating either from a geometric or from a combinatorial perspective.

Modern Schubert calculus seeks to construct the coefficients $c^{\nu}_{\lambda, \mu}$ explicitly for various different groups $G$, parabolic subgroups $P$, or cohomology theories.  One of the most vexing open questions in the field is the one closest to our starting point: to find an explicit, positive construction of the coefficients in the cohomology ring $H^*(GL_n(\mathbb{C}/B))$ of the flag variety of type $A$.

\subsection{Classical representation theory of Springer fibers}\label{subsection: Springer representations} The combinatorial parametrizations of the cells in Springer varieties hint at the representations that arise in the Springer context.  We open this section with a variation of Theorem \ref{theorem: all Springer cells} that makes this hint more explicit.  Recall that a standard Young tableau is one that is filled so that rows increase left-to-right and columns increase top-to-bottom.

\begin{theorem}[Spaltenstein 1976]
The top-dimensional cells $C_w \cap \mathcal{S}_X$ are bijective with the set of standard fillings of $\lambda(X)$.
\end{theorem}

Standard Young tableaux of shape $\lambda$ also count a classical representation-theoretic object: the dimension of the irreducible representation of $S_n$ associated the partition $\lambda$ of $n$.  This is not coincidence.

\begin{theorem}\label{theorem: Springer representation type A}
Fix a nilpotent matrix $X: \mathbb{C}^n \rightarrow \mathbb{C}^n$.  The cohomology $H^*(\mathcal{S}_X)$ of the Springer fiber carries an $S_n$-action.  The top-dimensional cohomology is irreducible and in fact is the irreducible representation associated to $\lambda(X)$.  Each irreducible representation of $S_n$ can be obtained uniquely in this way by varying over the conjugacy classes of nilpotent matrices.
\end{theorem}

Originally proven by Springer \cite{Spr76}, the theorem has many different proofs using different approaches and perspectives, some recovering only parts of the theorem as we have stated it and others much stronger.  For instance there are proofs due to Kazhdan and Lusztig \cite{KazLus80}, Borho and MacPherson \cite{BorMac83}, Lusztig \cite{Lus84}, Garsia and Procesi \cite{GarPro92}, and many others.  Of course there is a second irreducible $S_n$-representation of the same dimension as the irreducible representation associated to $\lambda$, namely its dual (obtained by tensoring with the sign representation).  Interestingly the literature on Springer's representation is ambiguous on this point: different constructions of ``the" Springer representation use either $\lambda$ or its dual.  Hotta appears to be the first to recognize this subtlety and classify different constructions up to that moment \cite{Hot81}. The representation described in this theorem generalizes to all Lie types, though Lusztig showed that outside of type $A$ the top-dimensional cohomology need not be bijective with irreducible representations; he defined {\em cuspidal representations} to be those that do not appear.  

A proof of Theorem \ref{theorem: Springer representation type A} is outside the scope of this survey but we include a sketch of Grothendieck's approach (see Grinberg's exposition for more \cite{Gri98}).  We consider two subspaces of the Lie algebra $\mathfrak{g}$: 
\begin{itemize}
\item the nilpotent subalgebra $\mathcal{N}$ consisting of all nilpotent elements of $\mathfrak{g}$ and \item the subalgebra $\mathfrak{g}^{rs}$ consisting of all regular semisimple elements of $\mathfrak{g}$.  
\end{itemize}
Now define the subspace $\widetilde{\mathfrak{g}}$ of the product space $\mathfrak{g} \times G/B$ by
\[\widetilde{\mathfrak{g}} = \{ (X, gB): g^{-1}Xg \in \mathfrak{b}\}\]
and define subspaces $\widetilde{\mathcal{N}}$ and $\widetilde{\mathfrak{g}^{rs}}$ analogously.  By projecting to the first factor we get the following commutative diagram.
\[\begin{array}{ccccc}
\widetilde{\mathcal{N}} & \hookrightarrow & \widetilde{\mathfrak{g}} & \hookleftarrow & \widetilde{\mathfrak{g}^{rs}} \\
\downarrow & & \downarrow & & \downarrow \\
{\mathcal{N}} & \hookrightarrow & {\mathfrak{g}} & \hookleftarrow & {\mathfrak{g}}^{rs}
\end{array}\]
This diagram has many useful and surprising features.  First the projection $\widetilde{\mathcal{N}} \rightarrow {\mathcal{N}}$ is a resolution of singularities, which we can see by projecting $\widetilde{\mathcal{N}}$ to the second factor and noting that the fiber over each flag $gB$ is isomorphic to $\mathfrak{b}$.  Moreover the fiber of the map $\widetilde{\mathcal{N}} \rightarrow {\mathcal{N}}$ over the element $X$ is the Springer fiber $\mathcal{S}_X$.   (This explains the name Springer {\em fiber}.)  On the right-hand side, the projection $\widetilde{\mathfrak{g}^{rs}} \rightarrow {\mathfrak{g}^{rs}}$ is an $n!$-sheeted cover on which the group $S_n$ acts as deck transformations. (The regular semisimple case in Example \ref{example: Springers} gives a small example of this cover.)  

More subtlely, the map $\widetilde{\mathcal{N}} \rightarrow {\mathcal{N}}$ has a geometric property called {\em semismall}, which is a condition that constrains the size of the fibers over the singular part of $\mathcal{N}$.  It allows us to use the decomposition theorem of Beilinson-Bernstein-Deligne-Gabber \cite{BBD82}, which in an informal sense breaks the total cohomology into pieces incorporating geometric subspaces from $\mathfrak{g}$ paired with $S_n$-actions on other pieces of $\mathfrak{g}$.  From a dimension count, the only pieces that survive are the geometric subspaces corresponding to conjugacy classes in $\mathcal{N}$ together with $S_n$-representations from the generic part of the Lie algebra.  The generic element of $\mathfrak{g}$ is regular semisimple, so the representations come from $\mathfrak{g}^{rs}$ which we observed carries the regular representation and thus decomposes as desired.  We emphasize that this proof relies on being able to calculate dimensions and other geometric properties of Springer fibers in order to prove hypotheses of the theorem, including that the map is semismall.  

We end this section with a remark about one way that this representation generalizes.  Hessenberg varieties are a larger family of varieties than Springer fibers, in which $\mathfrak{b}$ is replaced by different subspaces $H$ of the Lie algebra (more general even than parabolic subalgebras).  The Weyl group acts on the cohomology of nilpotent Hessenberg varieties by the {\em monodromy} representation, though the structure of the monodromy representation is more mysterious.  Shareshian-Wachs conjectured that the representation arises in the combinatorics of certain quasisymmetric functions that they studied for independent reasons \cite{ShaWac12}; recently both Brosnan and Chow \cite{BroCho} and Guay-Paquet \cite{Gua} proved this conjecture independently and with different methods. Another open question asks for explicit combinatorial descriptions of the representations that arise for various $X$ and $H$. 

\subsection{Non-classical representation theory and components of Springer fibers}\label{subsection: Khovanov Springers}  The previous two sections described geometric representations arising from the total space of a variety, namely in $H^*(G(k,n),\mathbb{C})$ in the Schubert case and on $H^*(\mathcal{S}_X)$ in the Springer case.  In this section we describe interesting new work relating certain quantum representations to the {\em components} of $\mathcal{S}_X$ in the type $A$ case. This section is less detailed than the previous two, partly because this representation was more recently discovered and our understanding of it continues to evolve.   

Spiders are diagrammatic categories encoding representations; the spider for $A_n$ encodes the representations of $\mathcal{U}_q(\mathfrak{sl_n})$.  The objects in this category represent tensors of representations.  The morphisms are planar graphs with boundary called {\em webs}, and the spider is equipped with skein-theoretic braiding morphisms that allow us to interpret tangles as webs.   In some descriptions the diagrams are drawn like tangles: the top and bottom represent the weight of the source and target representation and the twists in the strands indicate an action on the corresponding tensor products. 

\begin{figure}[h]
\begin{picture}(100,20)(0,0)
\put(-5,0){\curve(0,0, 10,10, 20,0)}
\put(15,0){\curve(10,0, 20,10, 30,0)}
\put(0,0){\curve(-15,0, 20,20, 55,0)}
\put(30,0){\curve(30,0, 40,10, 50,0)}

\put(-20,-10){\small{1}}
\put(-7,-10){\small{2}}
\put(12,-10){\small{3}}
\put(20,-10){\small{4}}
\put(42,-10){\small{5}}
\put(51,-10){\small{6}}
\put(59,-10){\small{7}}
\put(78,-10){\small{8}}
\end{picture}
\hspace{0.25in}
$\begin{array}{|c|c|c|c|}
\cline{1-4} 1 & 2 & 4 & 7 \\
\cline{1-4} 3 & 5 & 6 & 8 \\
\hline \end{array}$
\caption{A web when $\lambda(X)$ has two rows and the corresponding Young tableau} \label{figure: 2-row diagram}
\end{figure}
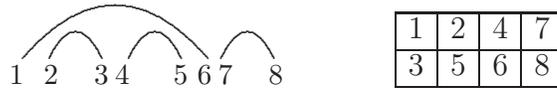

For instance in the $A_1$ case the webs are just Temperley-Lieb diagrams or crossingless matchings. When the partition $\lambda(X)$ has two rows, namely $X$ has two Jordan blocks, the components of the Springer fiber $\mathcal{S}_X$ naturally index a basis for the webs of $A_1$ spiders.  ``Naturally" here means that the combinatorics and the geometry of the components reflect information encoded in the diagram of the corresponding morphism. More precisely Khovanov proved that each component of the Springer fiber $\mathcal{S}_X$ is homeomorphic to an iterated tower of $\mathbb{P}^1$-bundles, where the structure of the tower corresponds to the nesting of arcs in the web corresponding to the component \cite{Kho04}.   Figure \ref{figure: 2-row diagram} gives an example of web for $A_1$ and its corresponding tableau; the corresponding Springer component is the product of  $\mathbb{P}^1$ with a $\mathbb{P}^1$-bundle over $\mathbb{P}^1 \times \mathbb{P}^1$.  Khovanov  discovered this connection between Springer fibers and webs when he was analyzing a ring that arose in his construction of Khovanov homology for tangles; he happened to observe that the center of his ring was isomorphic to the cohomology $H^*(\mathcal{S}_X)$ using work of Fung \cite{Fun03} and Garsia-Procesi \cite{GarPro92}.  

\begin{figure}[h]
\begin{picture}(80,20)(0,0)
\put(-10,0){\curve(0,0, 15,15, 30,0)}
\put(20,0){\curve(10,0, 25,15, 40,0)}

\put(5,0){\line(0,1){15}}
\put(45,0){\line(0,1){15}}

\put(-12,-10){\small{1}}
\put(3,-10){\small{2}}
\put(18,-10){\small{3}}
\put(27,-10){\small{4}}
\put(42,-10){\small{5}}
\put(57,-10){\small{6}}
\end{picture}
$\begin{array}{|c|c|}
\cline{1-2} 1 & 4  \\
\cline{1-2} 2 & 5  \\
\cline{1-2}   3 &  6   \\
\hline \end{array}$
\hspace{1in}
\begin{picture}(80,20)(15,0)
\put(0,0){\curve(0,0, 10,10, 20,0)}
\put(5,0){\curve(25,0, 35,10, 45,0)}
\put(10,0){\curve(50,0, 60,10, 70,0)}

\put(10,10){\curve(0,0, 30,10, 60,0)}
\put(40,10){\line(0,1){10}}

\put(-2,-10){\small{1}}
\put(18,-10){\small{2}}
\put(28,-10){\small{3}}
\put(48,-10){\small{4}}
\put(58,-10){\small{5}}
\put(78,-10){\small{6}}
\end{picture}
$\begin{array}{|c|c|}
\cline{1-2} 1 & 3  \\
\cline{1-2} 2 & 5  \\
\cline{1-2}   4 &  6   \\
\hline \end{array}$
\caption{Webs when $\lambda(X)$ has three rows and their corresponding Young tableaux} \label{figure: 3-row diagram}
\end{figure}
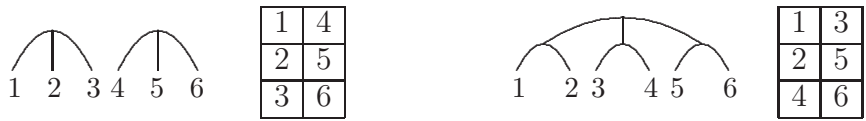

When $\lambda(X)$ has three rows Khovanov and Kuperberg showed that the standard Young tableaux combinatorially index the morphisms for $A_2$ spiders \cite{KhoKup99}.   It seems that there is a geometric relationship between the combinatorics of the diagram and the components of the corresponding Springer fiber: for instance, a diagram without any crossings corresponds to a smooth component (and in fact an iterated fiber product of copies of $GL_3(\mathbb{C})/B)$).  The example on the left in Figure \ref{figure: 3-row diagram} shows a component that is homeomorphic to $GL_3(\mathbb{C})/B \times GL_3(\mathbb{C})/B$ while the example on the right is a more complicated component.  One open question asks for a deeper analysis of how combinatorial features in the webs  (crossings, nesting, etc.) correspond to geometric properties of the components (singularities, nested product structures, etc.).  

We close this section by sketching a categorical description of the connection between Springer fibers and $A_n$ spiders, one which extends to a larger theory connecting similar varieties to other representation-theoretic constructions.  
In the case of spiders, we want a map between $(n,m)$-tangles and isomorphism classes of certain exact functors $D_n \rightarrow D_m$ that preserves knot-theoretic and algebraic structures of the tangles.  Cautis and Kamnitzer proved that if $D_n$ is the derived category of certain equivariant coherent sheaves on Springer fibers of shape $(n,n)$ then braid moves and other tangle operations correspond to Fourier-Mukai transforms \cite{CauKam08}.  Part of their insight is that Springer fibers naturally arise in the geometric Langlands program, and in fact are associated via the geometric Satake correspondence to the same $\mathfrak{sl}_n$ representations that arise from the knot-theoretic perspective; this is the core idea that they generalize to other varieties. 

\section{Connecting Springer fibers with Schubert varieties}\label{section: Springers and Schuberts}
The previous sections have used the analogy between Springer fibers and Schubert varieties to describe their geometry and representation theory in more depth.  In this section, we describe a web of new results and conjectures that connects the two kinds of varieties more deeply and directly.

To begin we define permutations $w_T$ associated to standard Young tableaux $T$.  

\begin{definition}
Fix a standard Young tableau $T$ with $n$ boxes.  For each $i$ with $2 \leq i \leq n$ 
\begin{itemize}
\item let $d_i$ denote the number of rows strictly above $i$ in $T$ and
\item let $w_i$ denote the increasing product of simple transpositions
\[s_{i-d_i} s_{i-d_i+1} \cdots s_{i-2} s_{i-1}\]
where each $s_i = (i,i+1)$.  (Our convention is that if $i=0$ then $w_i=e$ is the identity.)
\end{itemize}
Then the Schubert point associated to $T$ is the permutation 
\[w_T = w_n w_{n-1} w_{n-2} \cdots w_2\]
\end{definition}

\begin{example}\label{example: finding w_T}
For instance we have the following:
\[\begin{array}{|c|c|} \cline{1-2} 1 & 4  \\ \cline{1-2} 2 & 5  \\ \cline{1-2} 3 & \multicolumn{1}{c}{} \\ \cline{1-1} \end{array} \hspace{0.5in} \longleftrightarrow \hspace{0.5in} \begin{array}{l} d_5 = 1 \\ d_4 = 0 \\ d_3 = 2 \\ d_2 = 1 \\ d_1 = 0 \end{array} \hspace{0.5in} \longleftrightarrow \hspace{0.5in} \left(s_4\right)(s_1s_2)(s_1)\]
\end{example}

It turns out that these permutations $w_T$ index a set of Schubert varieties whose union has the same Betti numbers as a Springer fiber.  More precisely we have the following \cite{PreTym1}.

\begin{theorem}[Tymoczko-Precup] \label{theorem: TymPre}
The Betti numbers of the Springer variety associated to $X$ agree with the Betti numbers of the union of Schubert varieties
\[H_*\left(\mathcal{S}_X\right) \cong H_*\left(\bigcup \overline{C_{w_T}}\right)\]
where the union is taken over all standard tableaux $T$ of shape $\lambda(X)$.
\end{theorem}

\begin{example}
To show the subtlety of this result, we will calculate the Euler characteristic of each side of this equation in the case when $\lambda = (2,2,1)$.  In other words we will count the number of nonempty cells $C_w \cap \mathcal{S}_X$ and the number of nonempty cells in the union $\bigcup \overline{C_{w_T}}$ not worrying about degrees except when it is unavoidable. 

To find the Euler characteristic of the Springer fiber, we need to count the number of row-strict tableaux of shape $(2,2,1)$ by Theorem \ref{theorem: all Springer cells}.  The bottom box can be filled with any number, so there are five choices for that entry.  The four numbers left must then fill the shape $(2,2)$ so that each row increases. In other words we just partition the four numbers into pairs and order each pair in the unique increasing way.  There are $\binom{4}{2}=6$ ways to divide four distinct numbers into pairs and so $5 \cdot 6=30$ row-strict tableaux overall.

To get the Euler characteristic of the union of Schubert varieties, we need to compute all standard tableaux of shape $(2,2,1)$.  There are five in this case, which we found using brute force and the fact that $1$ must go in the top-left box while $5$ must go in one of the two bottom-right corners.
\[\begin{array}{|c|c|} \cline{1-2} 1 & 4  \\ \cline{1-2} 2 & 5  \\ \cline{1-2} 3 & \multicolumn{1}{c}{} \\ \cline{1-1} \end{array}  \hspace{0.5in}
\begin{array}{|c|c|} \cline{1-2} 1 & 3  \\ \cline{1-2} 2 & 5  \\ \cline{1-2} 4 & \multicolumn{1}{c}{} \\ \cline{1-1} \end{array}  \hspace{0.5in}
\begin{array}{|c|c|} \cline{1-2} 1 & 2  \\ \cline{1-2} 3 & 5  \\ \cline{1-2} 4 & \multicolumn{1}{c}{} \\ \cline{1-1} \end{array}  \hspace{0.5in}
\begin{array}{|c|c|} \cline{1-2} 1 & 3  \\ \cline{1-2} 2 & 4  \\ \cline{1-2} 5 & \multicolumn{1}{c}{} \\ \cline{1-1} \end{array}  \hspace{0.5in}
\begin{array}{|c|c|} \cline{1-2} 1 & 2  \\ \cline{1-2} 3 & 4  \\ \cline{1-2} 5 & \multicolumn{1}{c}{} \\ \cline{1-1} \end{array} 
\]

We created $w_T$ for the first of these tableaux in Example \ref{example: finding w_T}.  All of the Schubert points $w_T$ are found in the same way.  We list them below in the same order as their corresponding tableaux:
\[s_4s_1s_2s_1, s_4s_2s_3s_1, s_4s_2s_3s_2, s_3s_4s_3s_1, s_3s_4s_3s_2\]
To enumerate the cells in a particular Schubert variety, we just need to determine all the subwords of each permutation.  The challenge of determining the cells in a union of Schubert varieties is that many of these subwords can coincide. For instance the Schubert variety corresponding to the permutation $s_3s_4s_3s_2$ is the union of the Schubert cells corresponding to
\[e, s_2, s_3, s_4, s_3s_4, s_4s_3, s_2s_4, s_3s_2, s_3s_4s_3, s_3s_4s_2, s_4s_3s_2, s_3s_4s_3s_2\]
(omitting duplicates and using relations in the permutation group to simplify where possible).  But all of those except $s_3s_4s_3s_2$ and $s_3s_4s_2$ index a Schubert cell in one of the Schubert varieties  earlier in our list.  In other words, the Schubert variety for $s_3s_4s_3s_2$ only increases the Euler characteristic of the union by two.  (In this calculation, the dimension of each cell comes ``for free" since it's simply the number of simple reflections in the word.)  Performing this calculation for  the whole union, we get the following Betti numbers:
\[1, 4, 9, 11, 5\]
Their sum is 30, consistent with the theorem.
\end{example}

This result is one piece of a more general collection of conjectures and results involving Hessenberg varieties, the same generalization of Springer fibers mentioned in Section \ref{subsection: Springer representations}.  Indeed it appears that the parameters of a nilpotent $X$ and a Hessenberg space $H$ determine a union of Schubert varieties whose homology is the same as the corresponding Hessenberg variety. Harada and the author proved this conjecture for {\em Peterson varieties}, when $X$ is regular nilpotent (meaning has a single Jordan block) and $H$ consists of the subspace of matrices that are zero below the subdiagonal \cite{HarTym11}.  Mbirika proved the conjecture when $X$ is regular nilpotent and $H$ is arbitrary \cite{Mbi10}.  The case of the Springer fiber is more complicated and is also the a key tool in a more recent result by the author and Precup, which extends it to the case when $H$ is parabolic (the so-called parabolic Hessenberg varieties) \cite{PreTym2}.  Some of these results also extend to general Lie type \cite{HarTym11}.

We want to stress that this result only describes an enumerative and combinatorial property of Springer fibers.  It says nothing about, for instance, the cohomology class induced by $\mathcal{S}_X$ in $H^*(G/B)$ or the multiplicative structure of the ring $H^*(\mathcal{S}_X)$, though those are both interesting questions.  Indeed this conjecture came about in part because of work to determine the equivariant cohomology ring of Hessenberg varieties.  In particular, it appears that the equivariant cohomology of the Hessenberg variety corresponding to $X$ and $H$ can be determined by restricting a set of equivariant Schubert classes $\mathcal{Y}$ to certain fixed points in the Hessenberg variety.  The equivariant cohomology of the specific union of Schubert varieties identified by these conjectures appears to be computed by restricting the same set of equivariant Schubert classes $\mathcal{Y}$ instead to fixed points in the Schubert varieties.  That is an exciting conjecture though attempts to prove it have so far been limited to special cases like that of Peterson varieties \cite{HarTym11} by algebraic and combinatorial challenges.   

For these reasons we conjecture that an underlying geometric principle determines these results, even though all known proofs are combinatorial.  Indeed we conjecture that there is a degeneration of nilpotent Hessenberg varieties to unions of Schubert varieties, perhaps similar to Knutson-Miller's degeneration of Schubert varieties to unions of line bundles \cite{KnuMil05}.

\section{Open questions}

For convenience we list here the open questions throughout the whole paper.

\begin{itemize}
   \item Characterize the intersections $C_w \cap \mathcal{S}_X$ for arbitrary nilpotent $X$.  For which $X$ in a fixed conjugacy class is $C_w \cap \mathcal{S}_X$ affine?  For which $X$ are the nonempty intersections $C_w \cap \mathcal{S}_X$ enumerated by row-strict tableaux of shape $\lambda(X)$?  (See Section \ref{section: geometry of Springer cells}.)
    \item Give a complete list of the singular cells of Springer fibers, namely for each nilpotent $X$ and permutation flag $wB \in \mathcal{S}_X$, give a closed condition to determine if the closure $\overline{C_w \cap \mathcal{S}_X}$ is singular or not.  (See Section \ref{section: geometry of Springer cells}.)
    \item Determine which fixed points are in the closure of each Springer cell, namely for each nilpotent $X$ and pair of permutation flags $vB, wB \in \mathcal{S}_X$ give a closed condition to determine if $vB$ is in the closure of the cell $\overline{C_w \cap \mathcal{S}_X}$.  (See Section \ref{section: geometry of Springer cells}.)
    \item Determine which Springer Schubert cell closures $\overline{C_w \cap \mathcal{S}_X}$ are Cohen-Macaulay (Gorenstein, etc.).  (See Section \ref{section: geometry of Springer cells}.)
    \item Find an explicit, positive construction of the coefficients in the cohomology ring $H^*(GL_n(\mathbb{C}/B)$ of the flag variety of type $A$.  (See Section \ref{subsection: Schubert calculus}.)
    \item Find an explicit combinatorial description of the representation of $S_n$ on the cohomology of Hessenberg varieties for different $X$ and $H$. (See Section \ref{subsection: Springer representations}.)
    \item Identify the combinatorial features in web diagrams (crossings, nesting, etc.) that correspond to geometric properties in the closures of the corresponding components of the Springer fiber (singularities, nested product structures, etc.).  (See Section \ref{subsection: Khovanov Springers}.)
     \item Extend Theorem \ref{theorem: TymPre} to all nilpotent Hessenberg varieties.  (See Section \ref{section: Springers and Schuberts}.)
     \item Identify closed combinatorial formulas for the cohomology class induced by the Springer fiber $\mathcal{S}_X$ in $H^*(G/B)$.  (See Section \ref{section: Springers and Schuberts}.)
    \item Identify the structure constants of the cohomology $H^*(\mathcal{S}_X)$ of the Springer fiber in terms of the basis of Springer Schubert classes. (See Garsia-Procesi's work on $H^*(\mathcal{S}_X)$ \cite{GarPro92} as well as Section \ref{section: Springers and Schuberts}.)
     \item Identify the equivariant cohomology $H^*_T(\mathcal{S}_X)$ of the Springer fiber according to the outline in \cite{HarTym}, or any other way.  (See Section \ref{section: Springers and Schuberts}.)
     \item Give a geometric explanation for Theorem \ref{theorem: TymPre}.  (See Section \ref{section: Springers and Schuberts}.)
\end{itemize}

\bibliographystyle{plain}

\def\cprime{$'$}

\end{document}